# CARNE–VAROPOULOS BOUNDS FOR CENTERED RANDOM WALKS

By Pierre Mathieu

*CMI*

We extend the Carne–Varopoulos upper bound on the probability transitions of a Markov chain to a certain class of nonreversible processes by introducing the definition of a "centering measure." In the case of random walks on a group, we study the connections between different notions of centering.

**1. Introduction.** Let $X = (X_t, t \in \mathbb{N})$ be a Markov chain taking its values in some discrete set, $V$.

The paper is concerned with two related issues: in Section 2, the state space of the Markov chain is not assumed to have any special algebraic structure. We introduce a "centering condition" which generalizes the classical reversibility assumption. The main result is an extension of the Carne–Varopoulos inequality for the transition probabilities of a not necessarily reversible Markov chain; see Theorems 2.8 and 2.10. In Section 3 we restrict our attention to random walks on groups. We then investigate the relation between different possible definitions of a "centered random walk."

The initial motivation of this work was to find a different, more geometrical and combinatorial interpretation of the bounds obtained by Alexopoulos for random walks on nilpotent groups; see [1]. This is partially achieved, as far as the upper bound is concerned, in Proposition 3.3(a). But it turned out that our notion of *centering measure* can also be used to study nonreversible random walks on other examples of groups, such as Baumslag Solitar groups or wreath products; see Section 3.

*The Carne–Varopoulos bound.* A measure, $\pi$, on $V$ is called *reversible* for the Markov chain $X$ if the following detailed balance condition is satisfied:

for all $x, y \in V$,

(1) $$\pi(x)\mathbb{P}[X_1 = y|X_0 = x] = \pi(y)\mathbb{P}[X_1 = x|X_0 = y].$$

Not all Markov chains admit a reversible measure.

The detailed balance condition is equivalent to saying that the transition operator of $X$ is symmetric in $L^2(V, \pi)$. It is then possible to apply different tools from analysis, in particular spectral theory, to study the Markov chain. As an example of a distinguished property of reversible Markov chains, let us quote the Carne–Varopoulos upper bound: assume that $\pi$ is a reversible measure for $X$; then, for all $x, y \in V$ and $t \in \mathbb{N}^*$, we have

(2) $$\mathbb{P}[X_t = y|X_0 = x] \leq 2\sqrt{\frac{\pi(y)}{\pi(x)}} e^{-d^2(x,y)/(2t)}.$$

In (2), $d(x,y)$ is the natural distance associated to $X$, that is, the minimal number of steps required for the Markov chain to go from $x$ to $y$. The first paper to deal with such long-range estimates for transition probabilities is [11]. We refer to [3] or [13], Theorem 14.12 and Lemma 14.21 for a proof of (2) which relies on spectral theory. Inequality (2) gives a crude upper bound on the tail of the law of $X_t$ which turned out to be very useful in the analysis of the long-time behavior of reversible Markov chains.

*Centered random walks on a graph.* This paper arose as an attempt to get a similar bound for a not necessarily reversible Markov chain. Thus we do not assume that $X$ admits a reversible measure and ask: does there exist a constant $C$ such that, for all $x, y \in V$ and $t \in \mathbb{N}^*$,

(3) $$\mathbb{P}[X_t = y|X_0 = x] \leq Ce^{-d^2(x,y)/(Ct)}?$$

In the case of random walks in $\mathbb{Z}^d$, that is, if $X_t$ is obtained as a sum of $t$ independent, identically distributed random variables with finite support in $\mathbb{Z}^d$, then inequality (3) holds if and only if the mean value of $X_1$ vanishes or, equivalently, $E[X_t] = 0$ for all $t \in \mathbb{N}$. By analogy, we interpret (3) as a centering condition for the Markov chain $X$ although, for a general set $V$, it does not make sense anymore to speak of "vanishing mean" for $X_1$.

The transition probabilities of $X$ endow its state space $V$ with a structure of weighted oriented graph. In the second part of the paper, we define the class of *centered Markov chains* in terms of a splitting on this graph into oriented cycles; see Definition 2.1. Markov chains admitting a reversible measure are centered. We then prove a Carne–Varopoulos upper bound of the form (3) in Theorem 2.8. We also prove that the Dirichlet form satisfies a sector condition and derive some easy consequences in terms of Green kernels; see Lemma 2.12 and Proposition 2.13. In order to illustrate our definition, a special case of our general result is described at the end of this introduction.

*Centered random walks on a group.* The third part of the paper is devoted to random walks on groups. That is, we assume that $V$ is a discrete group; choose a finite generating set for $V$, say $G$ and define $X_t$ as a sum of independent, uniformly distributed random variables on $G$. Let $\mu$ be the uniform probability distribution on $G$, and let $\mu^t$ denote the $t$th convolution power of $\mu$. Thus $\mu^t$ is the law of $X_t$. In this context, (3) reads: does there exist a constant $C$ such that, for all $x \in V$ and $t \in \mathbb{N}^*$,

$$\mu^t(x) \leq Ce^{-d^2(id,y)/(Ct)}? \tag{4}$$

Here $id$ is the unit element in $V$. $d(x,y)$ is the word distance between $x$ and $y$. Up to multiplicative constants, $d(x,y)$ is independent of the choice of the generating set.

The graph associated to the random walk $X$ is now a Cayley graph of $V$, but, unless $G$ is symmetric, this is an oriented Cayley graph. Finding cycles in this Cayley graph amounts to writing $id$ as a product of elements of $G$. We may apply results of the second part to derive sufficient conditions on $G$ that imply (4): let $\mathcal{N}$ be the semigroup made of the elements of $V$ that can be written as products of elements in $G$ where each of the elements of $G$ appears the same number of times. In Proposition 3.1, we show that if $id \in \mathcal{N}$, then (4) is satisfied for some constant $C$. One can also consider sums of independent, identically distributed random variables with a more general law than the uniform distribution over $G$.

Checking whether $id \in \mathcal{N}$ is an—apparently new—combinatorial problem involving the geometry of $V$ and the choice of $G$. We solve it for nilpotent groups. Baumslag–Solitar groups, examples of wreath products and free groups are also considered; see Section 3.3.

As a consequence, in the above mentioned examples, we obtain the equivalence of the following two centering conditions:

(C1) $id \in \mathcal{N}$;
(C2) the image of the uniform measure on $G$ by any homomorphism of $V$ on $\mathbb{R}$ has vanishing mean.

*Application to the rate of escape.* Carne–Varopoulos bounds can be used in order to bound the rate of escape of the random walk from its initial point. In the case of a centered Markov chain, it is easy to deduce from the Carne–Varopoulos bound that the rate of escape vanishes if the volume growth is subexponential; see Theorem 2.11. In the case of random walks on a group, one can do much better and prove that the speed vanishes if and only if the Poisson boundary is trivial; see Proposition 3.11. This last statement extends well-known results for symmetric random walks; see [5, 9] or [12], among other references.

*An example.* We consider the special case of a Markov chain associated to an oriented unweighted graph structure on $V$. So let $E \subset V \times V$ be such that, for all $x \in V$, the number of points $y \in V$ such that $(x, y) \in E$ is finite and uniformly bounded in $x$. The Markov process $(X_t, t \in \mathbb{N})$ is defined by the usual rule: at each step, one selects at random (with uniform distribution) one of the edges in $E$ starting from the current position. Then the random walker jumps along the chosen edge.

A *cycle* is a sequence $\gamma = (x_0, x_1, \ldots, x_k)$ in $V$ such that $x_k = x_0$ and $(x_i, x_{i+1}) \in E$ for all $i = 0, \ldots, (k-1)$. We allow cycles of the form $(x_0, x_0)$ or $(x_0, x_1, x_0)$. Let $|\gamma| = k$ be the length of $\gamma$. We write that the edge $(x, y)$ belongs to $\gamma$ if, for some $i$, we have $x = x_i$ and $y = x_{i+1}$.

Assume that there exists a collection of cycles, $(\gamma_i, i \in \mathbb{N})$, satisfying the following two properties: (i) $\sup_i |\gamma_i| < \infty$, (ii) any edge $(x, y) \in E$ belongs to exactly one of the $\gamma_i$'s; then (3) holds for some constant $C$.

Now suppose that $V$ is a group with generating set $G = (g_1, \ldots, g_K)$. Then $E = \{(x, y) : x^{-1}y \in G\}$ defines an oriented Cayley graph on $V$. Cycles correspond to relations in $V$. Conditions (i) and (ii) are satisfied if there is a permutation of $\{1, \ldots, K\}$, say $\sigma$, such that $g_{\sigma(1)} \cdot g_{\sigma(2)} \cdots g_{\sigma(K)} = id$. Then (4) is satisfied.

The condition $g_{\sigma(1)} \cdot g_{\sigma(2)} \cdots g_{\sigma(K)} = id$ obviously implies that, for any homomorphism $h$ of $V$ on $\mathbb{R}$, $\sum_{i=1}^{K} h(g_i) = 0$. Whether the converse is true or not depends on the group; see Section 3.

*Further references.* The idea of using a decomposition of the state space of a Markov chain into cycles is not new. We refer in particular to the work of Kalpazidou [10] and to the first chapters of the book [7]. However, these authors are mostly interested in recurrent Markov chains.

The main technical tools used to prove our main result, Theorem 2.8, are borrowed from the work of Hebisch and Saloff-Coste, although some extra work is necessary to handle the lack of reversibility.

Comparison theorems for Green kernels similar to our Proposition 2.13(i) have been obtained by various authors; see, for instance, [2] or [4].

## 2. Centered Markov chains on graphs.

2.1. *Definitions.* In this section we introduce the definitions related to the graph structure induced by a Markov chain on its state space. As in the Introduction, let $(X_t, t \in \mathbb{N})$ be a Markov chain taking its values in some infinite countable set, $V$. We assume that $X$ is irreducible.

For $x$ and $y$ in $V$, define $q(x, y) = \mathbb{P}[X_1 = y | X_0 = x]$. Considering $q(x, y)$ as the weight of the edge $(x, y) \in V \times V$, we can see $\Gamma = (V, q)$ as a weighted, oriented graph.

Call a *cycle* a finite sequence $\gamma = (x_0, x_1, \ldots, x_k)$ of points in $V$ such that $x_k = x_0$ and $q(x_i, x_{i+1}) > 0$ for all $i = 0, \ldots, (k-1)$. We allow cycles of the form $(x_0, x_0)$ or $(x_0, x_1, x_0)$. Sometimes we identify the cycle $\gamma$ with a sequence of edges, that is, $\gamma = ((x_0, x_1), \ldots, (x_{k-1}, x_k))$. Define $|\gamma| = k$ to be the *length* of $\gamma$. We further suppose that cycles are edge self-avoiding, that is, that $(x_i, x_{i+1}) = (x_j, x_{j+1})$ implies that $i = j$. But we do not assume that cycles are vertex self-avoiding.

DEFINITION 2.1. Let $m$ be a measure on $V$. The graph $\Gamma$ is *centered* if there is a collection of cycles $(\gamma_i, i \in \mathbb{N})$ and positive weights $(q_i, i \in \mathbb{N})$ such that:

(i) $\sup_i |\gamma_i| < \infty$,
(ii) for any $x, y \in V$, we have

(5) $$m(x)q(x,y) = \sum_i q_i \mathbf{1}_{(x,y) \in \gamma_i}.$$

We then call $m$ a *centering measure* for the process $(X_t)$ (or for the graph $\Gamma$).

To avoid empty statements, we shall always assume that $m$ is not identically vanishing. From Remark 2.6 below it will follow that $m(x) > 0$ for all $x \in V$.

We shall use the notation $\varepsilon = \inf_{x \in V} m(x) \geq 0$ and $C_0 = \sup_i |\gamma_i|$.

REMARK 2.2. We may suppress the condition that cycles have to be edge self-avoiding. Let us call "*generalized cycle*" a sequence satisfying all the properties of cycles except it may have edge self-intersections. For a given edge, $(x, y) \in V \times V$, let $N((x, y), \gamma) = \#\{e \in \gamma : (x, y) = e\}$ be the number of occurrences of $(x, y)$ in the generalized cycle $\gamma$.

$\Gamma$ is then centered iff there exists a collection of generalized cycles, $(\gamma_i, i \in \mathbb{N})$, such that $\sup_i |\gamma_i| < \infty$ and, for all $x, y \in V$, we have

(6) $$m(x)q(x,y) = \sum_i q_i N((x,y), \gamma_i).$$

This fact is easy to prove by splitting generalized cycles into edge self-avoiding cycles.

REMARK 2.3 (The reversible case). Suppose that $m$ is a reversible measure for $X$, that is, assume that the detailed balance condition is satisfied: for any $x, y \in V$,

$$m(x)q(x,y) = m(y)q(y,x).$$

Choose cycles of the form $\gamma = (x, y, x)$ whenever $q(x, y) > 0$ and $\gamma = (x, x)$ whenever $q(x, x) > 0$. To the cycle $(x, y, x)$, we attach the weight $q = m(x)q(x, y)$; to the cycle $(x, x)$, we attach the weight $q = m(x)q(x, x)$. It is then immediate to deduce from the detailed balance condition that condition (5) holds. In other words, reversible graphs are centered.

EXAMPLE 2.4 (Unweighted graphs). *Let $E \subset V \times V$. Assume that, for all $y \in V$, the number of points $x \in V$ such that $(x, y) \in E$ is finite. Let $N^+(x) = \{y \in V : (x, y) \in E\}$, and define*

$$q(x, y) = \begin{cases} \dfrac{1}{\#N^+(x)}, & \text{if } (x, y) \in E, \\ 0, & \text{otherwise,} \end{cases}$$

*so that the random walker moves by choosing uniformly at random an edge in $E$ starting from its current position and then jumping along the chosen edge. Let $m(x) = \#N^+(x)$.*

*Assume that there exists a collection of cycles, $(\gamma_i, i \in \mathbb{N})$, and an integer, $n$, such that* (i) $\sup_i |\gamma_i| < \infty$, *and* (ii) *for any edge $e \in E$, $\#\{i : e \in \gamma_i\} = n$. Then $\Gamma$ is centered.*

PROOF. Indeed we have

$$\sum_i \mathbf{1}_{(x,y) \in \gamma_i} = n = nm(x)q(x, y),$$

for any edge $(x, y) \in E$. Thus we may choose the weights $q_i = \frac{1}{n}$ to check condition (5). □

Note that, for $\Gamma$ to be centered for the measure $m$, it is necessary that $\#\{y \in V : (y, x) \in E\} = \#\{y \in V : (x, y) \in E\}$ for all $x \in V$.

LEMMA 2.5. *Let $\Gamma$ be centered for $m$. Then $m$ is an invariant measure for $X$, that is, for all $y \in V$, one has $\sum_{x \in V} m(x)q(x, y) = m(y)$.*

PROOF. For given $x \in V$ and $i \in \mathbb{N}$, note that there exists $y \in V$ with $(x, y) \in \gamma_i$ iff there exists $y \in V$ with $(y, x) \in \gamma_i$. Because cycles are edge self-avoiding, $\#\{y \in V : (x, y) \in \gamma_i\} = \#\{y \in V : (y, x) \in \gamma_i\}$. Therefore

$$\sum_y \sum_i q_i \mathbf{1}_{(x,y) \in \gamma_i} = \sum_y \sum_i q_i \mathbf{1}_{(y,x) \in \gamma_i}.$$

Thus

$$\sum_x m(x)q(x, y) = \sum_x \sum_i q_i \mathbf{1}_{(x,y) \in \gamma_i}$$

$$= \sum_x \sum_i q_i \mathbf{1}_{(y,x) \in \gamma_i}$$
$$= \sum_x m(y)q(y,x) = m(y). \qquad \square$$

REMARK 2.6. As a consequence of the lemma, since we have assumed that $X$ is irreducible, we must have $m(x) > 0$ for all $x \in V$. Keeping in mind that the weights $q_i$ are positive, we note that it implies that, for any $x, y \in V$, $q(x, y) > 0$ if and only if there exists at least one $i \in \mathbb{N}$ such that $(x, y) \in \gamma_i$.

We now recall the definition of the distance associated to $\Gamma$. For $x, y \in V$, let $d(x, y)$ be the smallest $k \in \mathbb{N}$ such that there is a sequence $x_0, \ldots, x_k$ with $x_0 = x$, $x_k = y$ and $q(x_i, x_{i+1}) + q(x_{i+1}, x_i) > 0$. In other words, $d$ is the classical graph distance associated to the undirected graph structure on $V$ defined by

$$E^0 = \{(x, y) \in V \times V : q(x, y) + q(y, x) > 0\}.$$

REMARK 2.7. Assume that $\Gamma$ is centered. If $d(x, y) = k$, then there exists a sequence $(x_0, \ldots, x_K)$ such that $x_0 = x$, $x_K = y$ and $q(x_i, x_{i+1}) > 0$, for all $i$. Besides we may choose $K \leq C_0 k$.

Indeed, if $d(x, y) = 1$, then, either $q(x, y) > 0$—and then $K = 1$—or $q(x, y) = 0$, in which case $q(y, x) > 0$. In the latter case, we choose one cycle $\gamma_i$ such that $(y, x) \in \gamma_i$, say $\gamma_i = (y, x, x_2, \ldots, x_{a-1}, y)$. Then $a \leq C_0$. Besides we have found a path, $(x, x_2, \ldots, x_{a-1}, y)$, of length bounded by $a \leq C_0$, linking $x$ to $y$ and such that $q(e_1, e_2) > 0$ when $(e_1, e_2) \in \gamma_i$. Thus the claim is proved for $k = 1$. The general case follows.

We can now state the main result of this section:

THEOREM 2.8. *Let $\Gamma$ be a centered graph for the measure $m$. Assume that $\varepsilon = \inf_{x \in V} m(x) > 0$. Then there exists a constant $C$, that only depends on $\varepsilon$ and $C_0$, such that, for all $x, y \in V$ and $t \in \mathbb{N}^*$, we have*

$$\mathbb{P}[X_t = y | X_0 = x] \leq Cm(y)e^{-d^2(x,y)/(Ct)}.$$

2.2. *Proof of Theorem* 2.8.

*Preliminaries on Dirichlet forms.* Define the operator $Qf(x) = \mathbb{E}[f(X_1) | X_0 = x] = \sum_{y \in V} q(x, y)f(y)$ on functions with finite support. $Q_t$ will denote the $t$th power of $Q$.

Let $Q^*$ be the adjoint of $Q$ with respect to the measure $m$. Then $Q^*f(x) = \sum_{y \in V} q^*(x,y)f(y)$, with $q^*(y,x) = \frac{m(x)}{m(y)}q(x,y)$. Using (5), we get that

$$q^*(y,x)m(y) = \sum_i q_i \mathbf{1}_{(x,y) \in \gamma_i}.$$

This last formula may as well be written

$$q^*(x,y)m(x) = \sum_i q_i \mathbf{1}_{(x,y) \in \gamma_i^*},$$

where, for a cycle $\gamma$, we use the notation $\gamma^*$ to denote the reversed cycle. (Reverse the order of the sequence defining $\gamma$.) Thus the graph $\Gamma^* = (V, q^*)$ is also centered for the same measure $m$. It is actually the graph associated to the time reversal of the Markov chain $X$. In particular all the results we are about to prove for centered graphs may be applied to $\Gamma^*$.

We have already noticed that $m(Qf) = m(f)$. The operator $Q$ being positivity preserving, we thus have $m(|Qf|) \leq m(|f|)$. It is also clear that $\sup_{x \in V} |Qf(x)| \leq \sup_{x \in V} |f(x)|$. It follows from Jensen's inequality, or by interpolation, that $Q$ is a contraction in $L^p(V,m)$ for all $p \in [1, \infty]$. By duality, $Q^*$ is also a contraction in $L^p(V,m)$.

Define the Dirichlet form $\mathcal{E}(f,g) = m(g.(I-Q)f)$. It can be expressed with the kernel $q$ by

$$\mathcal{E}(f,g) = \sum_{x,y \in V} m(x)q(x,y)g(x)(f(x)-f(y)).$$

We also consider the symmetrized Dirichlet form

$$\mathcal{E}^0(f,g) = \frac{1}{2}(\mathcal{E}(f,g) + \mathcal{E}(g,f))$$
$$= m\left(g.\left(I - \frac{Q+Q^*}{2}\right)f\right)$$
$$= \sum_{x,y \in V} m(x)\frac{q(x,y)+q^*(x,y)}{2}g(x)(f(x)-f(y)).$$

Since, $m(x)(q(x,y) + q^*(x,y)) = m(x)q(x,y) + m(y)q(y,x)$, we have

(7) $\qquad \mathcal{E}^0(f,g) = \frac{1}{2} \sum_{x,y \in V} p^0(x,y)(f(x)-f(y))(g(x)-g(y)),$

with

(8) $\qquad p^0(x,y) = p^0(y,x) = \frac{1}{2}(m(x)q(x,y) + m(y)q(y,x)).$

Let us now compute the antisymmetric part of $\mathcal{E}$:

$$\mathcal{E}(f,g) - \mathcal{E}^0(f,g) = m\left(g . \frac{Q^* - Q}{2} f\right)$$

$$= \sum_{x,y \in V} m(x) g(x) f(y) \frac{q^*(x,y) - q(x,y)}{2}$$

$$= \frac{1}{2} \sum_{x,y \in V} (f(x)g(y) - f(y)g(x)) m(x) q(x,y).$$

And, using (5), we obtain the useful representation formula:

(9) $$\mathcal{E}(f,g) - \mathcal{E}^0(f,g) = \tfrac{1}{2} \sum_i q_i \sum_{(x,y) \in \gamma_i} (f(x)g(y) - f(y)g(x)).$$

*Poincaré inequality.* We shall use the following Poincaré inequality on the discrete circle: let $\gamma$ be a cycle. There exists a constant, $C_\gamma$, such that, for all functions $g$ such that $\sum_{x \in \gamma} g(x) = 0$, we have

(10) $$\sum_{x \in \gamma} g(x)^2 \leq C_\gamma \sum_{(x,y) \in \gamma} (g(x) - g(y))^2.$$

The best constant in (10) is the inverse spectral gap of the nearest-neighbor symmetric random walk on $\gamma$; thus (10) is a Poincaré inequality. Besides, the constant $C_\gamma$ depends only on the length $|\gamma|$.

PROOF OF THEOREM 2.8. The "symmetric" version of Theorem 2.8 is stated as Theorem 14.12 in [13]. (The argument is due to Hebish and Saloff-Coste; see [6].) We try to follow Woess as closely as possible, starting with the next lemma, but there is an extra nonsymmetric term to be handled by specific arguments. This is where the assumption (5) enters into play.

Keep in mind that $C$ is a constant which is allowed to depend only on $\varepsilon$ and $C_0$. Choose some reference point $o \in V$. For $s \in \mathbb{R}$, define the function $w_s(x) = e^{sd(o,x)}$. We need the following.

LEMMA 2.9. *There exists a constant $C$, that depends on $C_0$ only, and such that, for all $s \in \mathbb{R}$, $|s| \leq \frac{1}{C}$, and for any function $f$ with finite support, we have*

$$\mathcal{E}(w_s f, w_{-s} f) \geq -Cs^2 (1 + e^{C|s|}) m(f^2).$$

PROOF. We use the notation $w = w_s$ and note that, replacing $f$ with $wf$, we have to prove that

$$\mathcal{E}(w^2 f, f) \geq -Cs^2(1 + e^{C|s|}) m(w^2 f^2).$$

Using expression (7), we get that

$$4\mathcal{E}(w^2 f, f) = 4\mathcal{E}^0(w^2 f, f) + 4\mathcal{E}(w^2 f, f) - 4\mathcal{E}^0(w^2 f, f)$$
$$= 2 \sum_{(x,y) \in V} p^0(x,y)(w(y)^2 f(y) - w(x)^2 f(x))(f(y) - f(x))$$
$$+ 4(\mathcal{E}(w^2 f, f) - \mathcal{E}^0(w^2 f, f))$$
$$= A_1 + A_2 + B,$$

where

$$A_1 = \sum_{(x,y) \in V} p^0(x,y)(f(x) - f(y))^2 (w(x)^2 + w(y)^2),$$
$$A_2 = \sum_{(x,y) \in V} p^0(x,y)(f(x)^2 - f(y)^2)(w(x)^2 - w(y)^2),$$
$$B = 4(\mathcal{E}(w^2 f, f) - \mathcal{E}^0(w^2 f, f)).$$

From the proof of Lemma 14.14 in [13], we have

(11) $$(A_2)^2 \le 8s^2(1 + e^{2|s|}) A_1 m(w^2 f^2).$$

We need a similar estimate for $B$. We first rewrite $B$ using the set of paths $(\gamma_i, i \in \mathbb{N})$ as in (9):

$$B = 2 \sum_i q_i \sum_{(x,y) \in \gamma_i} f(x) f(y) (w(x)^2 - w(y)^2).$$

For $i \in \mathbb{N}$, we use the notation $c_i$ for the mean value of $f$ on the points of the cycle $\gamma_i$, and $f_i(x) = f(x) - c_i$. Taking into account that $\gamma_i$ is a closed path shows that $\sum_{(x,y) \in \gamma_i} w(x)^2 - w(y)^2 = 0$. Therefore

$$B = 2 \sum_i q_i \bigg( \sum_{(x,y) \in \gamma_i} f_i(x) f_i(y) (w(x)^2 - w(y)^2)$$
$$+ 2c_i \sum_{(x,y) \in \gamma_i} (f_i(x) + f_i(y))(w(x)^2 - w(y)^2) \bigg).$$

If $d(x,y) = 1$, then $|w(x) - w(y)| \le C|s|(w(x) + w(y))$. Writing $\underline{w}_i$ (resp. $\overline{w}_i$) for the min (resp. max) of $w$ over the path $\gamma_i$, we have

$$|B| \le C|s| \sum_i q_i (\overline{w}_i)^2 \bigg( \bigg( \sum_{x \in \gamma_i} |f_i(x)| \bigg)^2 + |c_i| \sum_{x \in \gamma_i} |f_i(x)| \bigg).$$

We now use the Poincaré inequality (10) for the function $f_i$ to deduce that

$$\left(\sum_{x\in\gamma_i}|f_i(x)|\right)^2 \leq |\gamma_i|\sum_{x\in\gamma_i}(f_i(x))^2 \leq C_{\gamma_i}|\gamma_i|\sum_{(x,y)\in\gamma_i}(f_i(x)-f_i(y))^2$$
$$= C_{\gamma_i}|\gamma_i|\sum_{(x,y)\in\gamma_i}(f(x)-f(y))^2.$$

The length of $\gamma_i$ being bounded by $C_0$, we may therefore choose a constant $C$, independent of $i$, such that

$$\left(\sum_{x\in\gamma_i}|f_i(x)|\right)^2 \leq C\sum_{(x,y)\in\gamma_i}(f(x)-f(y))^2.$$

Also note that $(c_i)^2 \leq C\sum_{x\in\gamma_i}f^2(x)$.

From the previous inequalities, we conclude that

$$|B| \leq C|s|\sum_i q_i(\overline{w}_i)^2\Bigg(\sum_{(x,y)\in\gamma_i}(f(x)-f(y))^2$$
$$+ \sqrt{\sum_{x\in\gamma_i}f^2(x)}\sqrt{\sum_{(x,y)\in\gamma_i}(f(x)-f(y))^2}\Bigg).$$

For the next step, we use the fact that $w$ is roughly constant on each path $\gamma_i$. More precisely, since $|\gamma_i| \leq C_0$, two points on $\gamma_i$ are at distance at most $C_0$. Therefore $\overline{w}_i \leq e^{C|s|}\underline{w}_i$, where $C$ depends only on $C_0$. Therefore

$$|B| \leq C|s|e^{C|s|}\sum_i q_i\Bigg(\sum_{(x,y)\in\gamma_i}(f(x)-f(y))^2(w(x)^2+w(y)^2)$$
$$+ \sqrt{\sum_{x\in\gamma_i}f^2(x)w(x)^2}$$
$$\times \sqrt{\sum_{(x,y)\in\gamma_i}(f(x)-f(y))^2(w(x)^2+w(y)^2)}\Bigg)$$
$$\leq C|s|e^{C|s|}\Bigg(\sum_i q_i \sum_{(x,y)\in\gamma_i}(f(x)-f(y))^2(w(x)^2+w(y)^2)$$
$$+ \sqrt{\sum_i q_i \sum_{x\in\gamma_i}f^2(x)w(x)^2}$$
$$\times \sqrt{\sum_i q_i \sum_{(x,y)\in\gamma_i}(f(x)-f(y))^2(w(x)^2+w(y)^2)}\Bigg),$$

where we used the Cauchy–Schwarz inequality. Using (5), we deduce that

$$|B| \leq C|s|e^{C|s|}\bigg(\sum_{x,y \in V}(f(x)-f(y))^2(w(x)^2+w(y)^2)m(x)q(x,y)$$
$$+\sqrt{\sum_{x \in V}f(x)^2w(x)^2\sum_i q_i \mathbf{1}_{x \in \gamma_i}}$$
$$\times \sqrt{\sum_{x,y \in V}(f(x)-f(y))^2(w(x)^2+w(y)^2)m(x)q(x,y)}\bigg).$$

But $m(x) = \sum_{y \in V} m(x)q(x,y) = \sum_i q_i \sum_{y \in V} \mathbf{1}_{(x,y) \in \gamma_i} \geq \sum_i q_i \mathbf{1}_{x \in \gamma_i}$ and $m(x)q(x,y) \leq 2p^0(x,y)$. Therefore

$$|B| \leq C|s|e^{C|s|}\bigg(\sum_{x,y \in V}(f(x)-f(y))^2(w(x)^2+w(y)^2)p^0(x,y)$$
$$+\sqrt{\sum_{x \in V}f(x)^2w(x)^2m(x)}$$
$$\times \sqrt{\sum_{x,y \in V}(f(x)-f(y))^2(w(x)^2+w(y)^2)p^0(x,y)}\bigg),$$

that is,

(12) $$|B| \leq C|s|e^{C|s|}(A_1 + \sqrt{A_1 m(w^2 f^2)}).$$

Inequalities (11) and (12) clearly imply the lemma. □

We shall not explain how to deduce the theorem from the lemma since the arguments can be copied from the proof of Theorem 14.12 in [13]. (A referee pointed out that this is true up to the following additional observation: in the middle of page 156 of [13] one reads: "the adjoint of $P_s$ is $P_{-s}$." This is not the case here but everything applies to $Q^*$ in place of $Q$.) As in Theorem 14.12 in [13], we have in fact proved the stronger result:

THEOREM 2.10. *Let $\Gamma$ be a centered graph for the measure $m$. Assume that $\varepsilon = \inf_{x \in V} m(x) > 0$. Assume that there are constants $C_1$ and $d \geq 0$ such that, for all $x, y \in V$ and all $t \in \mathbb{N}^*$, we have*

(13) $$\mathbb{P}[X_t = y | X_0 = x] \leq C_1 m(y) t^{-d/2}.$$

*Then there exists a constant $C$ that only depends on $\varepsilon$, $d$, $C_0$ and $C_1$, such that, for all $x, y \in V$ and $t \in \mathbb{N}^*$, we have*

$$\mathbb{P}[X_t = y | X_0 = x] \leq Cm(y)t^{-d/2}e^{-d^2(x,y)/(Ct)}.$$

Theorem 2.8 is only the special case of Theorem 2.10 when $d = 0$. □

2.3. *Rate of escape.* The next statement is an easy consequence of the Carne–Varopoulos bounds.

THEOREM 2.11. *Assume that $\Gamma$ is centered for a measure $m$ such that $\varepsilon = \inf_{x \in V} m(x) > 0$. Let $V(t) = \sharp\{x \in V : d(o,x) \leq t\}$ be the volume of the ball centered at $o$. ($o$ is an arbitrary reference point.) If $\limsup_{t \to \infty} \frac{1}{t} \log V(t) = 0$, then for all $\alpha > 0$, we have*

$$\lim_{t \to +\infty} \mathbb{P}[d(o, X_t) \geq \alpha t | X_0 = o] = 0.$$

PROOF. Use Theorem 2.8 and the fact that $d(o, X_t) \leq t$ if $X_0 = o$, to get that

$$\mathbb{P}[d(o, X_t) \geq \alpha t | X_0 = o] = \sum_{x; \alpha t \leq d(o,x) \leq t} \mathbb{P}[X_t = x | X_0 = o]$$

$$\leq \sum_{x; \alpha t \leq d(o,x) \leq t} C e^{-d^2(o,x)/(Ct)}$$

$$\leq C e^{-\alpha^2 t / C} V(t) \to 0. \qquad \square$$

2.4. *Sector condition and Green kernels.*

LEMMA 2.12 (Sector condition). *Let $\Gamma$ be a centered graph for the measure $m$. There exists a constant $M$, function of $C_0$ only, such that, for all finitely supported functions $f$ and $g$, we have*

$$\mathcal{E}(f,g)^2 \leq M^2 \mathcal{E}(f,f) \mathcal{E}(g,g).$$

PROOF. We write that $\mathcal{E}(f,g) = \mathcal{E}^0(f,g) + \mathcal{E}(f,g) - \mathcal{E}^0(f,g)$, where, as before, $\mathcal{E}^0(f,g) = \frac{1}{2}(\mathcal{E}(f,g) + \mathcal{E}(g,f))$ is the symmetric part of $\mathcal{E}$.

Since $\mathcal{E}^0$ is a symmetric bilinear form, we have $\mathcal{E}^0(f,g)^2 \leq \mathcal{E}^0(f,f)\mathcal{E}^0(g,g) = \mathcal{E}(f,f)\mathcal{E}(g,g)$. It remains to prove that $(\mathcal{E}(f,g) - \mathcal{E}^0(f,g))^2 \leq M^2 \mathcal{E}(f,f)\mathcal{E}(g,g)$.

From (9), we know that

$$\mathcal{E}(f,g) - \mathcal{E}^0(f,g) = \frac{1}{2} \sum_i q_i \sum_{(x,y) \in \gamma_i} (f(x)g(y) - f(y)g(x)).$$

Note that the quantity $\sum_{(x,y) \in \gamma_i}(f(x)g(y) - f(y)g(x))$ remains unchanged if we modify by a constant the value of $f$ or $g$ on $\gamma_i$. Thus let $c_i$ (resp. $d_i$) be the mean of $f$ (resp. $g$) on $\gamma_i$ and set $f_i = f - c_i$ (resp. $g_i = g - d_i$). From the Poincaré inequality (10), we get a constant $M_i$, that depends on the length of $\gamma_i$ only, such that

$$\sum_{x \in \gamma_i} f_i^2(x) \leq M_i \sum_{(x,y) \in \gamma_i} (f(x) - f(y))^2,$$

$$\sum_{x \in \gamma_i} g_i^2(x) \leq M_i \sum_{(x,y) \in \gamma_i} (g(x) - g(y))^2.$$

Since the length of $\gamma_i$ is bounded by $C_0$, we have $M = \sup_i M_i < \infty$. Then

$$\left(\sum_{(x,y)\in\gamma_i} (f(x)g(y) - f(y)g(x))\right)^2$$

$$= \left(\sum_{(x,y)\in\gamma_i} (f_i(x)g_i(y) - f_i(y)g_i(x))\right)^2$$

$$\leq M^2 \sum_{x\in\gamma_i} f_i^2(x) \sum_{x\in\gamma_i} g_i^2(x)$$

$$\leq M^2 \sum_{(x,y)\in\gamma_i} (f(x) - f(y))^2 \sum_{(x,y)\in\gamma_i} (g(x) - g(y))^2,$$

and therefore

$$(\mathcal{E}(f,g) - \mathcal{E}^0(f,g))^2$$

$$\leq M^2 \left(\sum_i q_i \sum_{(x,y)\in\gamma_i} (f(x)-f(y))^2\right)\left(\sum_i q_i \sum_{(x,y)\in\gamma_i} (g(x)-g(y))^2\right).$$

It now only remains to note that $\sum_i q_i \sum_{(x,y)\in\gamma_i}(f(x)-f(y))^2 = \sum_{x,y\in V}(f(x)-f(y))^2 q(x,y)m(x) = \mathcal{E}(f,f)$. □

We can use the sector condition of Lemma 2.12 to compare the Green kernel of the Markov chain $X$ with the Green kernel of the Markov chain associated to $\mathcal{E}^0$, say $X^0$. Let $Q^0 = \frac{Q+Q^*}{2}$. The operator $Q^0$ is then symmetric with respect to $m$ and has kernel $q^0(x,y) = \frac{1}{2}(q(x,y) + \frac{m(y)}{m(x)}q(y,x))$. By definition of the Dirichlet forms $\mathcal{E}$ and $\mathcal{E}^0$, one has the relation

$$m(f.(I-Q)f) = \mathcal{E}(f,f)$$
$$= m(f.(I-Q^0)f) = \mathcal{E}^0(f,f).$$

We use the notation $g(x,y)$ [resp. $g^0(x,y)$] to denote the Green kernel of $Q$ (resp. $Q^0$), be it finite or infinite. Thus

$$g(x,y) = \sum_{t\geq 0} \mathbb{P}[X_t = y | X_0 = x] = \frac{1}{m(x)} m(\delta_x.(I-Q)^{-1}\delta_y),$$

$$g^0(x,y) = \sum_{t\geq 0} \mathbb{P}[X_t^0 = y | X_0^0 = x] = \frac{1}{m(x)} m(\delta_x.(I-Q^0)^{-1}\delta_y).$$

PROPOSITION 2.13. (i) *For any* $x \in V$, *we have* $g(x,x) \leq g^0(x,x)$.

(ii) *Assume that $\Gamma$ is centered. Then, for all $x \in V$, $g^0(x,x) \leq M^2 g(x,x)$, where $M$ is the same constant as in Lemma* 2.12.

(iii) *As a consequence, if $\Gamma$ is centered, then $X$ is recurrent if and only if $X^0$ is recurrent.*

PROOF. Part (i) directly follows from Lemma 2.24 in [13] using the fact that $m(f.(I-Q)f) = m(f.(I-Q^0)f)$.

Part (ii) follows from Lemma (2.12):

$$\begin{aligned}
&(m(x)g^0(x,x))^2 \\
&= m((I-Q^0)^{-1}\delta_x.\delta_x)^2 \\
&= m((I-Q^0)^{-1}\delta_x.(I-Q)(I-Q)^{-1}\delta_x)^2 \\
&= \mathcal{E}((I-Q)^{-1}\delta_x, (I-Q^0)^{-1}\delta_x)^2 \\
&\leq M^2 \mathcal{E}((I-Q)^{-1}\delta_x, (I-Q)^{-1}\delta_x)\mathcal{E}((I-Q^0)^{-1}\delta_x, (I-Q^0)^{-1}\delta_x) \\
&= M^2 m(\delta_x.(I-Q)^{-1}\delta_x)m(\delta_x.(I-Q^0)^{-1}\delta_x) \\
&= M^2 m(x)g(x,x)m(x)g^0(x,x),
\end{aligned}$$

where we used Lemma 2.12 from line 4 to line 5. □

### 3. Centered Markov chains on groups.

3.1. *Definitions.* We shall apply the results of the previous section to the analysis of (nonreversible) random walks on groups. Our main purpose is to discuss the connections between different "natural" definitions of what a centered random walk on a group should be. Proposition 3.1 gives a simple sufficient condition for a random walk to be centered in the sense of Definition 2.1, and motivates the introduction of the centering condition (C1). We also consider the weaker but somehow more natural centering condition (C2). One question is then to decide whether, for a given group, conditions (C1) and (C2) are equivalent or not. We take up this problem in two steps: Section 3.2 contains some easy remarks on conditions (C1) and (C2) and a technical tool, Lemma 3.8, that turns out to be useful to deduce (C1) from (C2). In Section 3.3 we discuss different examples of groups. Finally, in Section 3.4 we prove that the velocity of a centered random walk vanishes if and only if its entropy also vanishes.

We therefore assume that $V$ is a discrete, infinite group of finite type and choose a finite sequence, $G = (g_1, \ldots, g_K)$ of elements of $V$. Note that we really mean a sequence, that is, the same element may appear more than once in $G$. *id* will denote the unit element in $V$. We say that $G$ is *generating* if the semigroup generated by $G$ is $V$: any element in $V$ can be written as a product of elements in $G$.

To $G$, we associate a Markov chain, $(X_t, t \in \mathbb{N})$, in the usual way: let $(U_i, i \in \mathbb{N}^*)$ be a sequence of independent random variables with uniform distribution in $\{1, \ldots, K\}$. Let $\eta_i = g_{U_i}$. We define the sequence $(X_t, t \in \mathbb{N})$ by the recursion relations:

$$X_0 = id,$$
$$X_{t+1} = X_t . \eta_{t+1}.$$

Let $\mathbb{P}$ be the law of the sequence $(X_t, t \in \mathbb{N})$. The law of $X_1$, say $\mu$, is easily computed:

$$\mu(x) = \frac{\#\{i : g_i = x\}}{K}.$$

The law of $X_t$ is then the $t$th convolution power of $\mu$, that we denote by $\mu^t$.

In the language of the first part of the paper, $X$ is the Markov chain associated to the graph $\Gamma = (V, q)$ with $q(x, y) = \frac{1}{K} \#\{i : g_i = x^{-1} \cdot y\}$.

We choose for reference measure $m$, the counting measure on $V$.

We recall that a function $\sigma : \{1, \ldots, nK\} \to \{1, \ldots, K\}$ is said to be *n to 1* if for all $i \in \{1, \ldots, K\}$, then $\#\{j \in \{1, \ldots, nK\} : \sigma(j) = i\} = n$.

PROPOSITION 3.1.  *We assume that there exist an integer $n \in \mathbb{N}^*$ and a function $\sigma : \{1, \ldots, nK\} \to \{1, \ldots, K\}$, which is n to 1 such that*

(14) $$g_{\sigma(1)} \cdot g_{\sigma(2)} \cdots g_{\sigma(nK)} = id.$$

*Then the graph $\Gamma$ is centered for the counting measure $m$. In particular, the conclusions of Theorem* 2.8 *and Lemma* 2.12 *hold.*

For further references, let us make a definition out of (14): we shall say that a given sequence $G$ satifies condition (C1) if there exist an integer $n \in \mathbb{N}^*$ and a function $\sigma : \{1, \ldots, nK\} \to \{1, \ldots, K\}$, which is $n$ to 1 and satisfies

(15) $$g_{\sigma(1)} \cdot g_{\sigma(2)} \cdots g_{\sigma(nK)} = id.$$

PROOF OF PROPOSITION 3.1.  Let $\tilde{g}_t = g_{\sigma(1)} \cdots g_{\sigma(t)}$ and let $\gamma_1$ be the cycle

$$\gamma_1 = ((id, \tilde{g}_1), (\tilde{g}_1, \tilde{g}_2), \ldots, (\tilde{g}_{nK-1}, \tilde{g}_{nK})).$$

By assumption $\tilde{g}_{nK} = id$. Also define the translated cycles: $\gamma_x = x.\gamma_1$, for all $x \in V$.

Because the cycles $\gamma_x$ may not be edge self-avoiding, we will use Remark 2.2 in Section 2.1 and check (6).

Let $a, b \in V$. The number of times the edge $(a, b)$ appears in a path $\gamma_x$ is the number of couples $(x, i)$ with $x \in V$ and $i \leq nK - 1$ and such that

$$(a, b) = (x.\tilde{g}_i, x.\tilde{g}_{i+1}),$$

or, equivalently,

(16) $$a = x.\tilde{g}_i \quad \text{and} \quad b = a.g_{\sigma(i+1)}.$$

If $q(a,b) = 0$, that is, $a_\cdot^{-1}b \notin G$, then (16) has no solution. Otherwise, $i$ being given, $x$ is uniquely determined by (16). Thus we are actually looking for the number of $i$'s such that $a_\cdot^{-1}b = g_{\sigma(i+1)}$. This number is $nKq(x,y)$, as clearly follows from the definition of $q$ and the assumption of $\sigma$ being $n$ to 1. □

We now introduce a second centering condition: a given sequence satisfies condition (C2) if for some integer $n$, $(g_1 \cdots g_K)^n \in [V,V]$ or, equivalently, $\sum_i h(g_i) = 0$ for any homomorphism $h$ from $V$ to $\mathbb{R}$; see Remark 3.6 below.

Note that the condition $(g_1 \cdots g_K)^n \in [V,V]$ is independent of the order in which the product is computed. Indeed, changing the order in this product would only multiply the result by an element in $[V,V]$.

Although condition (C1) is the one we needed to prove our results, condition (C2) is, to a certain extent, more natural. In particular, it is easier to check in examples.

It is also easy to see that (C1) implies (C2): indeed assume that (C1) holds. Then, since $\sigma$ is $n$ to 1, we obtain the product $(g_1 \cdots g_K)^n$ as a reordering of the elements of the product in (15). But changing the order in some product only multiplies this product by an element in $[V,V]$. Therefore $(g_1 \cdots g_K)^n \in [V,V]$ and (C2) holds.

DEFINITION 3.2. We will say that the group $V$ satisfies property (C) if, for any finite generating sequence, conditions (C1) and (C2) hold or fail simultaneously. *In extenso*, $V$ satisfies property (C) if, for any finite generating sequence $G = (g_1, \ldots, g_K)$ such that for some $n \in \mathbb{N}^*$ we have $(g_1 \cdots g_K)^n \in [V,V]$, then there exist an integer $n \in \mathbb{N}^*$ and a function $\sigma : \{1, \ldots, nK\} \to \{1, \ldots, K\}$, which is $n$ to 1 and satisfies

(17) $$g_{\sigma(1)} \cdot g_{\sigma(2)} \cdots g_{\sigma(nK)} = id.$$

PROPOSITION 3.3.

(a) *Nilpotent groups satisfy property* (C).
(b) *The Baumslag–Solitar group $BS_q$ satisfies property* (C).
(c) *The wreath product $\mathbb{Z} \wr \mathbb{Z}$ satisfies property* (C).
(d) *The free group $F_2$ does not satisfy property* (C).

REMARK 3.4. From Proposition 3.3(a) and property (C) it follows that if a generating set on a nilpotent group satisfies condition (C2), it then satisfies the Carne–Varopoulos upper bound. As a matter of fact, it would also be possible to use Theorem 2.10 to get an upper bound of the form $\mathbb{P}[X_t = y | X_0 = x] \leq Ct^{-r/2}e^{-d^2(x,y)/(Ct)}$ for any centered random walk. (Here $r$ is the volume growth exponent of the group.) But it should be pointed out that a more precise version of this last bound, and the corresponding lower bound, were obtained by Alexopoulos in [1] for more general centered random walks than ours. Alexopoulos' method is quite different from ours and does not use the equivalence between conditions (C1) and (C2).

3.2. *Centering conditions.* We start with some easy remarks on conditions (C1) and (C2):

REMARK 3.5. The random walk $X$, associated to the finite sequence $G$, lives on the semigroup generated by $G$. If (C1) holds, it is easy to see that the semigroup generated by $G$ is in fact a group.

REMARK 3.6 (Homomorphisms on $\mathbb{R}$). Let $G = (g_1, \ldots, g_K)$ be a finite sequence of elements of $V$.

First assume that for some $n$, $(g_1 \cdots g_K)^n \in [V, V]$. Then, for any homomorphism $h$ from $V$ to $\mathbb{R}$, we have $\sum_i h(g_i) = 0$.

Conversely, assume that, for any homomorphism from $V$ to $\mathbb{R}$, we have $\sum_i h(g_i) = 0$. Then, for some $n$, $(g_1 \cdots g_K)^n \in [V, V]$.

Indeed, let $\gamma$ be the image of the product $g_1 \cdots g_K$ on $V/[V, V]$. Either $\gamma$ has finite order—in which case the proof is finished—or it has infinite order. Since $V/[V, V]$ is Abelian, there exists a homomorphism $\tilde{h}$ from $V/[V, V]$ to $\mathbb{R}$ such that $\tilde{h}(\gamma) = 1$. Then $\tilde{h}$ induces a homomorphism on $V$ such that $h(g_1 \cdots g_K) = 1$. This is in contradiction with the assumption that $\sum_i h(g_i) = 0$.

Thus we have proved that, for a given sequence $G = (g_1, \ldots, g_K)$, the following two properties are equivalent:

(i) there exists $n$ such that $(g_1 \cdots g_K)^n \in [V, V]$,
(ii) for any homomorphism $h$ from $V$ to $\mathbb{R}$, we have $\sum_i h(g_i) = 0$.

REMARK 3.7. There are obvious counterexamples to the implication (C2) $\Longrightarrow$ (C1) for nongenerating sequences: choose $K = 1$. The condition (C1) is then equivalent to saying that $g_1$ has finite order. Condition (C2) is satisfied if $g_1 \in [V, V]$. Thus if $g_1 \in [V, V]$ but $g_1$ is of infinite order, then (C2) is satisfied but (C1) is not. We avoid this situation by assuming that the set $G$ generates $V$. Let us recall that the meaning of "generating" is: all elements of $V$ belong to the semigroup generated by $G$, that is, any $x \in V$ can be written as a product of elements in $G$.

The aim of the next section is to check that property (C) holds for some simple enough groups. The proofs are based on the following combinatorial lemma:

LEMMA 3.8. *Choose a finitely generated group, $V$, and some element $a \in V$. The following two properties are equivalent:*

(i) *For any finite generating sequence, $G = (g_1, \ldots, g_K)$, such that $(g_1 \cdots g_K)^n \in [V, V]$ for some $n \in \mathbb{N}^*$, then* (C1) *holds.*

(ii) *For any finite generating sequence, $G = (g_1, \ldots, g_K)$, such that $(g_1 \cdots g_K)^n \in [V, V]$ for some $n \in \mathbb{N}^*$, then* (C1) *holds for the enlarged sequence $(g_1, \ldots, g_K, a, a^{-1})$.*

PROOF. Of course (i) implies (ii). Assume that (ii) is verified. Let $G$ be some finite generating sequence such that $(g_1 \cdots g_K)^n \in [V, V]$. We check that $G$ satisfies (C1).

Since $G$ generates $V$, we can write
$$a = g_{\sigma_1(1)} \cdots g_{\sigma_1(k_1)} \quad \text{and} \quad a^{-1} = g_{\sigma_2(1)} \cdots g_{\sigma_2(k_2)},$$
for some applications $\sigma_1 : \{1, \ldots, k_1\} \to \{1, \ldots, K\}$ and $\sigma_2 : \{1, \ldots, k_2\} \to \{1, \ldots, K\}$. Call $G_1$ the sequence of elements in $V$ obtained by forming all the products of elements of $G$ of length $k_1$. In other words,
$$G_1 = (g_{\sigma_1(1)} \cdots g_{\sigma_1(k_1)}; \sigma_1 \in \{1, \ldots, K\}^{\{1,\ldots,k_1\}}).$$
(Remember $G_1$ is a sequence, not a set. The same element may appear more than once.) Similarly, define $G_2$ to be the sequence of elements in $V$ obtained by forming all the products of elements of $G$ of length $k_2$:
$$G_2 = (g_{\sigma_2(1)} \cdots g_{\sigma_2(k_2)}; \sigma_2 \in \{1, \ldots, K\}^{\{1,\ldots,k_2\}}).$$
Thus $a \in G_1$ and $a^{-1} \in G_2$. Finally let $\tilde{G}$ be the concatenation of the sequences $G$, $G_1$ and $G_2$. Then $\tilde{G}$ has $\tilde{K} = K + K^{k_1} + K^{k_2}$ elements.

We claim that $\tilde{G}$ satisfies the requirements of (ii). Indeed $\tilde{G}$ generates $V$ since it contains $G$ and $G$ generates $V$. We also have $a, a^{-1} \in \tilde{G}$. If we form the $n$th power of the product of the elements of $\tilde{G}$, we get
$$(\tilde{g}_1 \cdots \tilde{g}_{\tilde{K}})^n = ((g_1 \cdots g_K)\Pi_{\sigma_1}(g_{\sigma_1(1)} \cdots g_{\sigma_1(k_1)})\Pi_{\sigma_2}(g_{\sigma_2(1)} \cdots g_{\sigma_2(k_2)}))^n$$
$$= (g_1 \cdots g_K)^{n(1+k_1 K^{k_1-1}+k_2 K^{k_2-1})} \mod([V, V]),$$
where the second equality holds up to reordering.

Since, by assumption, $(g_1 \cdots g_K)^n \in [V, V]$, we see that $(\tilde{g}_1 \cdots \tilde{g}_2)^n \in [V, V]$. Therefore we deduce that $\tilde{G}$ satisfies the condition (C1): there exist some number $\tilde{n}$ and an application $\tilde{\sigma} : \{1, \ldots, \tilde{n}\tilde{K}\} \to \{1, \ldots, \tilde{K}\}$ such that

(18) $$\tilde{g}_{\tilde{\sigma}(1)} \cdots \tilde{g}_{\tilde{\sigma}(\tilde{n}\tilde{K})} = id,$$

and $\tilde{\sigma}$ is $\tilde{n}$ to 1. Imagine you rewrite the product (18) with the elements of $G$. From the construction of $\tilde{G}$, it then follows that each element of $G$ will appear exactly $\tilde{n}(1 + K^{k_1-1} + K^{k_2-1})$ times. We have thus checked condition (C1) for the generating sequence $G$.

Note that all over this proof the roles of the different elements of $G$ are symmetric. □

3.3. *Examples and proof of Proposition* 3.3.  As a preliminary, let us first consider the simplest example:

EXAMPLE 3.9 (Periodic groups).  We assume that all elements of $V$ have finite order. Thus $V$ is a periodic group, also called a torsion group. Given any finite set $G$, we can choose $n$ such that $g_i^n = id$ for all $i \in \{1, \ldots, K\}$. We then define $\sigma(i)$ to be the integer part of $1 + \frac{i-1}{n}$. $\sigma$ is clearly $n$ to 1. Besides,

$$g_{\sigma(1)} \cdots g_{\sigma(nK)} = g_1^n \cdots g_K^n = id.$$

As a conclusion the graph $\Gamma$ associated to $G$ is centered.

We shall now extend this result to more general random walks on $V$: let $\mu$ be a probability measure on $V$ with finite support. Consider the Markov chain with transition rates $q(x,y) = \mu(x_.^{-1}y)$. For $x \in V$ and $g$ in the support of $\mu$, define the cycle $\gamma_{x,g} = (x, x.g, x.g^2, \ldots, x.g^{p(g)})$, where $p(g)$ is the order of $g$. Let $q_g = \frac{1}{p(g)}\mu(g)$.

Choose $a, b \in V$. For fixed $g$, count the total number of occurrences of the edge $(a,b)$ in cycles of the form $\gamma_{x,g}$, where $x$ ranges through $V$. We get: $p(g)$ if $a_.^{-1}b = g$ and 0 otherwise. Therefore

$$\sum_{x,g} q_g N((a,b), \gamma_{x,g}) = \sum_g \mu(g) \mathbf{1}_{a_.^{-1}b=g} = \mu(a_.^{-1}b) = q(a,b).$$

We have checked condition (6) and therefore the graph $\Gamma = (V, q)$ is centered for the counting measure.

EXAMPLE 3.10 (Abelian case).  Assume that the product $g_1 \cdots g_K$ has finite order, say $p$. Then it is easy to construct a function $\sigma$, which is $p$ to 1 and satisfies (15). Now assume that $V$ is Abelian. If (C1) holds, then we must have $(g_1 \cdots g_K)^n = 1$. [This is just a reordering of the product in (15).] Then $g_1 \cdots g_K$ has finite order. Thus we see that, for an Abelian group, (C1) is fullfilled if and only if the product $g_1 \cdots g_K$ has finite order. In particular, Abelian groups satisfy property (C).

PROOF OF PROPOSITION 3.3.  (a) Nilpotent groups satisfy property (C). We proceed by induction on the nilpotency class of $V$. Let $V = V_0 > V_1 > \cdots > V_r = \{id\}$ be the lower central series of $V$ with $V_{i+1} = [V, V_i]$. Let $Z$ be the center of $V$. The case $r = 1$ corresponds to an Abelian group $V$ and was already discussed in Example 3.10.

Note that $V_{r-1}$ is Abelian and finitely generated. We may, and do, choose elements $(x_i, y_i, i = 1, \ldots, k)$ such that the set $([x_i, y_i], i = 1, \ldots, k)$ generates $V_{r-1}$. Finally notice that $V_{r-1} \subset Z$. Therefore if $x, y \in V$ are such that $[x, y] \in V_{r-1}$, then $[x^\alpha, y^\beta] = [x, y]^{\alpha\beta}$ for all nonnegative $\alpha$ and $\beta$.

Assume now that the statement of the proposition is true for any nilpotent group of class $r-1$ or less. Let $V$ be of class $r$. Let $G = (g_1, \ldots, g_K)$ be a finite generating sequence and let $n$ be such that $(g_1 \cdots g_K)^n \in [V, V]$. We wish to prove that condition (C1) holds. Using Lemma 3.8, it is sufficient to prove that the sequence $(g_1, \ldots, g_K, x_1, \ldots, x_k, y_1, \ldots, y_k, x_1^{-1}, \ldots, x_k^{-1}, y_1^{-1}, \ldots, y_k^{-1})$ satisfies (C1).

We use the induction assumption: the group $V/V_{r-1}$ is nilpotent of class strictly less than $r$. Therefore there is an integer $p$ and a 1 to $p$ function $\sigma : \{1, \ldots, pK\} \to \{1, \ldots, K\}$ such that $g_{\sigma(1)} \cdots g_{\sigma(pK)} \in V_{r-1}$. Therefore there exist $l \geq 0$, $l \leq k$ and $\alpha_1, \ldots, \alpha_l \in \mathbb{Z}$ such that $g_{\sigma(1)} \cdots g_{\sigma(pK)} [x_1, y_1]^{\alpha_1} \cdots [x_l, y_l]^{\alpha_l} = id$. Interchanging the roles of $x_i$ and $y_i$ when necessary, we may assume that the $\alpha_i$'s are nonnegative.

Let $\alpha$ be the product $\alpha = \alpha_1 \cdots \alpha_l$. Note that $[x_i, y_i]^{\alpha_i \alpha} = [x_i^{\alpha_i}, y_i^{\alpha_i}]^{\alpha/\alpha_i}$. We have

$$id = (g_{\sigma(1)} \cdots g_{\sigma(pK)} [x_1, y_1]^{\alpha_1} \cdots [x_l, y_l]^{\alpha_l})^{\alpha}$$
$$= (g_{\sigma(1)} \cdots g_{\sigma(pK)})^{\alpha} [x_1, y_1]^{\alpha_1 \alpha} \cdots [x_l, y_l]^{\alpha_l \alpha}$$
$$\text{(because } [x_i, y_i] \in V_{r-1} \subset Z\text{)}$$
$$= (g_{\sigma(1)} \cdots g_{\sigma(pK)})^{\alpha} [x_1^{\alpha_1}, y_1^{\alpha_1}]^{\alpha/\alpha_1} \cdots [x_l^{\alpha_l}, y_l^{\alpha_l}]^{\alpha/\alpha_l}$$
$$= (g_{\sigma(1)} \cdots g_{\sigma(pK)})^{\alpha} [x_1^{\alpha_1}, y_1^{\alpha_1}]^{\alpha/\alpha_1} \cdots [x_l^{\alpha_l}, y_l^{\alpha_l}]^{\alpha/\alpha_l} (x_1 x_1^{-1} y_1 y_1^{-1})^{\alpha(p-1)}$$
$$\times \cdots \times (x_l x_l^{-1} y_l y_l^{-1})^{\alpha(p-1)} (x_{l+1} x_{l+1}^{-1} y_{l+1} y_{l+1}^{-1})^{\alpha p} \cdots (x_k x_k^{-1} y_k y_k^{-1})^{\alpha p}.$$

In this last expression, each $g_i$ appears $\alpha p$ times; each term of the form $x_i$, $x_i^{-1}$, $y_i$ or $y_i^{-1}$ with $i \leq l$ appears $\alpha + \alpha(p-1) = \alpha p$ times; each term of the form $x_i$, $x_i^{-1}$, $y_i$ or $y_i^{-1}$ with $i > l$ appears $\alpha p$ times. Thus we have checked condition (C1).

(b) The Baumslag–Solitar group $BS_q$ satisfies property (C). By definition, the Baumslag–Solitar group $BS_q$ is the group with presentation $\langle a, b | ab = b^q a \rangle$, where $q \geq 2$ is an integer. It is an example of an amenable, solvable group of exponential volume growth. It is also the subgroup of the affine group of $\mathbb{R}$ generated by the transformations $x \to x + 1$ and $x \to qx$.

From the presentation, it is obvious that any homomorphism of $V$ on $\mathbb{R}$ should vanish on $b$. It is possible to prove that elements on $V$ can be written in the form $x = (a^{-l} b^m a^l) a^k$. In particular, if $x \in [V, V]$, then $x$ must be of the form $x = a^{-l} b^m a^l$ for some $l \geq 0$ and $m \in \mathbb{Z}$. We shall use the relation $[b^\beta, a^\alpha] = b^{\beta(1-q^\alpha)}$.

Let $G = (g_1, \ldots, g_K)$ be a finite generating sequence and choose $n$ such that $(g_1 \cdots g_K)^n \in [V, V]$. We wish to prove that condition (C1) holds. According to Lemma 3.8, it is sufficient to prove (C1) for the enlarged sequence $(g_1, \ldots, g_K, a, a^{-1}, b, b^{-1})$. Let $l \geq 0$ and $m \in \mathbb{Z}$ be such that $(g_1 \cdots g_K)^n = a^{-l} b^m a^l$.

First assume that $m \geq 0$. Choose $\alpha$ such that $q^\alpha - 1 \geq m$ and choose $j$ such that $j(q^\alpha - 1) \geq \alpha + l$. Let $k = j(q^\alpha - 1)$, $\beta = mj$, $k_1 = nk - \alpha - l \geq 0$ and $k_3 = nk - \beta \geq 0$.

We have
$$(g_1 \cdots g_K)^{kn} a^{-l} [b^\beta, a^\alpha] a^l a^{k_1} a^{-k_1} b^{k_3} b^{-k_3}$$
$$= (a^{-l} b^m a^l)^k a^{-l} b^{\beta(1-q^\alpha)} a^l = a^{-l} b^{km} b^{\beta(1-q^\alpha)} a^l = id,$$

since $km + \beta(1 - q^\alpha) = 0$.

Considering the expression $(g_1 \cdots g_K)^{kn} a^{-l} [b^\beta, a^\alpha] a^l a^{k_1} a^{-k_1} b^{k_3} b^{-k_3}$ as a word in the alphabet $G$, we see that: the elements $g_i, i \leq K$, appear each exactly $kn$ times; $a$ and $a^{-1}$ appear $\alpha + l + k_1 = kn$ times; $b$ and $b^{-1}$ appear $\beta + k_3 = kn$ times. Therefore we have checked (15).

The proof is done very much the same way if $m \leq 0$.

(c) The wreath product $\mathbb{Z} \wr \mathbb{Z}$ satisfies property (C). $\mathbb{Z} \wr \mathbb{Z}$ is isomorphic to the group of affine transformations of $\mathbb{R}$ generated by the translation $x \to x + 1$ and the homothety $x \to ax$ where $a$ is transcendental. It is also a semidirect product of $\mathbb{Z}$ and a direct product of countably many copies of $\mathbb{Z}$. It is therefore a two-step solvable group of finite type, although it is not finitely presented.

To be more precise, and quoting from [13]: a configuration $\eta$ is a function from $\mathbb{Z}$ to $\mathbb{Z}$ such that the set $\{x : \eta(x) \neq 0\}$ is finite. Equipped with pointwise addition, the set of configurations is a group, say $\tilde{\mathbb{Z}}$. $\mathbb{Z}$ acts on $\tilde{\mathbb{Z}}$ by automorphisms via $(y, \eta) \to T_y \eta$ where $T_y \eta(x) = \eta(x - y)$. The resulting semidirect product is the wreath product $\mathbb{Z} \wr \mathbb{Z}$. We denote by $\varepsilon$ the natural projection of $\mathbb{Z} \wr \mathbb{Z}$ onto $\mathbb{Z}$ and by $H$ the projection of $\mathbb{Z} \wr \mathbb{Z}$ on $\tilde{\mathbb{Z}}$. Thus any element of $\mathbb{Z} \wr \mathbb{Z}$ is a couple $a = (\varepsilon(a), \eta)$ where $\eta \in \tilde{\mathbb{Z}}$. $\mathbb{Z} \wr \mathbb{Z}$ is generated by the following four elements: $\tau_1 = (1, 0)$, $\tau_{-1} = \tau_1^{-1} = (-1, 0)$ and $\sigma_1 = (0, \eta_1)$, $\sigma_{-1} = \sigma_1^{-1} = (0, \eta_{-1})$, where $\eta_1(0) = 1$, $\eta_1(x) = 0$ if $x \neq 0$, $\eta_{-1}(0) = -1$, $\eta_{-1}(x) = 0$ if $x \neq 0$. We will use $|a|$ to denote the distance between $a \in \mathbb{Z} \wr \mathbb{Z}$ and $id$ in the metric induced by the generating set $\{\tau_1, \tau_{-1}, \sigma_1, \sigma_{-1}\}$.

$\varepsilon$ is a homomorphism of $\mathbb{Z} \wr \mathbb{Z}$ on $\mathbb{R}$. Another such homomorphism is $a \to \sum_{x \in \mathbb{Z}} H(a)(x)$.

Let $G$ be a finite generating set: $G = \{g_1, \ldots, g_K\}$. Assume that $G$ satisfies condition (C2). Therefore $\sum_{i=1}^K \varepsilon(g_i) = 0$ and $\sum_{i=1}^K \sum_{x \in \mathbb{Z}} H(g_i)(x) = 0$.

Using Lemma 3.8, in order to prove that condition (C1) holds we may, and will, replace $G$ by the enlarged generating set: $G' = \{\tau_1, \tau_1, \tau_{-1}, \tau_{-1}, \sigma_1, \sigma_{-1}, g_1, \ldots, g_K\}$.

We let $\phi$ be the product $\phi = g_1 \cdots g_K$ and $\phi_n = \phi(\tau_1.\phi)^n$. In the sequel to this proof, $C$ and $M$ will denote some constants that depend on $G$ but not on $n$.

We first note that $\varepsilon(\phi_n) = n$, since $\varepsilon(\phi) = 0$. Also note that $H(\phi_n)(x) = \sum_{j=0}^n H(\phi)(x - j) = \sum_{j=x-n}^x H(\phi)(j)$. And since $\sum_{x \in \mathbb{Z}} H(\phi)(x) = 0$, then there must be a constant $M$ such that $H(\phi_n)(x) \neq 0$ implies that $x \in [-M, M]$ or $x \in [n - M, n + M]$. Thus $\phi_n$ is of the form $\phi_n = A^{(n)} \tau_1^n B^{(n)}$ for some elements $A^{(n)}$ and $B^{(n)}$ such that $|A^{(n)}| + |B^{(n)}| \leq C$, for some constant $C$ (that does not depend on $n$!). Which means that we can write both $A^{(n)}$ and $B^{(n)}$ as products of elements of $\{\tau_1, \tau_{-1}, \sigma_1, \sigma_{-1}\}$ with less than $C$ symbols.

Thus we have obtained a trivial product: $id = \phi_n (B^{(n)})^{-1} \tau_{-1}^n (A^{(n)})^{-1}$ in which:

(i) each element $g_i$ appears $n+1$ times;
(ii) the numbers of occurrences of $\sigma_1$ and $\sigma_{-1}$ are equal because $\sum_{x \in \mathbb{Z}} H(\phi_n)(x) = 0 = \sum_{x \in \mathbb{Z}} H(A^{(n)} B^{(n)})(x)$. Call this number $b^{(n)}$. $b^{(n)}$ is bounded by some constant that does not depend on $n$, since $|A^{(n)}| + |B^{(n)}| \leq C$;
(iii) by the same argument, $\tau_1$ and $\tau_{-1}$ appear the same number of times, say $a^{(n)}$ and $a^{(n)} \leq n + C$.

Choose $n$ such that $b^{(n)} \leq n$ and $a^{(n)} \leq 2n + 2$. We obviously have $id = \phi_n B^{(n)-1} \tau_{-1}^n A^{(n)-1} (\tau_1 \tau_{-1})^{2n+2-a^{(n)}} (\sigma_1 \sigma_{-1})^{n+1-b^{(n)}}$ and this last expression proves (15).

(d) The free group $F_2$ does not satisfy property (C). Choose $G = (g_1 = a, g_2 = a^{-1}, g_3 = b, g_4 = b^{-1}, g_5 = b^{-2}, g_6 = ababa^{-2})$. Clearly, $G$ generates. Besides

$$g_1 g_6 g_3 g_5 g_2 g_4 = a^2 bab.a^{-2} b^{-1} a^{-1} b^{-1} = [a^2, bab] \in [V, V].$$

Let $n$ be a positive integer. Let $\gamma$ be an element of $V$ that can be written as a product of elements in $G$ using exactly $n$ times each of the $g_i$'s. Let us prove that $\gamma \neq id$.

First write $\gamma$ as a product of elements in $G$ with $n$ occurrences of each $g_i$. We label the different occurrences of $g_6$ by the numbers 1 to $n$ according to the order in which they appear. Replace the $g_i$'s by their expressions in terms of $a$, $b$, $a^{-1}$, $b^{-1}$. We obtain a nonreduced word in the alphabet $(a, b, a^{-1}, b^{-1})$. The letters coming from the $i$th occurrence of $g_6$ are labeled $i$. We run the following algorithm to reduce it step by step: read the word starting from the left; do all cancellations you find on your way; start again when you reach the end of the word. For $i, j \in \{1, \ldots, n\}$, we draw an edge between $i$ and $j$ if, while running the cancellation algorithm, one of the "$a^{-1}$" with label $i$ cancels with one of the "$a$" with label $j$ or one of the "$a^{-1}$" with label $j$ cancels with one of the "$a$" with label $i$. This way we obtain a nonoriented graph structure on $\{1, \ldots, n\}$. Let $J$ be the total number of edges of this graph. If $J < n$, then $\gamma$ is not $id$. Indeed, there are $3n$ occurrences of "$a^{-1}$" in the nonreduced word, $n$ of them coming from $g_2$ and $2n$ of them coming from $g_6$. Of the $2n$ occurrences of "$a^{-1}$" coming from $g_6$, $J$ cancel with some "$a$" coming from some occurrence of $g_6$, and, at most $n$ of them cancel with an "$a$" coming from $g_1$. Thus, after the algorithm has run, there will be at least $(2n-(n+J))$ "$a^{-1}$" left in the reduced word.

The graph structure we have built on $\{1,\ldots,n\}$ satisfies the following properties:

(i) it has no double edge, that is, we did not draw two edges from $i$ to $j$. This is due to the presence of the "$b$" between the two "$a$" in $g_6$;

(ii) it has no loop of the form $i \leftrightarrow j \leftrightarrow i$;

(iii) a configuration of the form $i_1 \leftrightarrow i_2$, $i_3 \leftrightarrow i_4$ with $i_4$ strictly between $i_1$ and $i_2$ implies that $i_3$ lies between $i_1$ and $i_2$ (in the broad sense).

(iii) follows from the definition of the algorithm.

Thus the graph contains no cycle. Indeed, if $i_1 \leftrightarrow i_2 \leftrightarrow \cdots \leftrightarrow i_k$ was a minimal cycle ($i_k = i_1$ and the labels $i_1,\ldots,i_{k-1}$ are pairwise different), then, from (iii), we deduce that, up to a circular permutation or running the cycle in the opposite order, the sequence $i_1,\ldots,i_{k-1}$ must be increasing. But this is impossible because the "$b$" would not cancel.

We conclude that the graph has no cycle. Therefore its number of edges is strictly less than $n$. $\square$

3.4. *On the velocity.* Given the finite generating set $G$, we consider the induced distance on $V$: $d(x,y)$ is the minimum number of elements in $G \cup G^{-1}$ whose product equals $x^{-1}y$. This definition corresponds to the definition of distance we used in Section 2.1.

The speed of the random walk $(X_t, t \in \mathbb{N})$ is $L = \lim_{t\to\infty} \frac{1}{t} d(id, X_t)$. The entropy of the random walk is $h = \lim_{t\to\infty} -\frac{1}{t} \log \mu^t(X_t)$ where $\mu$ is the law of $X_1$ (and therefore $\mu^t$ is the law of $X_t$). The subadditive ergodic theorem implies that the limits defining $L$ and $h$ exist in the almost sure sense as well as in the $L^1$ sense; both $L$ and $h$ are nonnegative numbers. (See [13], Theorem (8.14), [5], Section IV or [9], Theorem 1.6.4.)

It is known, without any symmetry assumption, that $h=0$ if and only if the Poisson boundary of the random walk is trivial; see [5], Section IV or [9], Section 1.6. From Corollary 1 in [12] it follows that $h = 0$ if $L = 0$. The converse follows from the classical Carne–Varopoulos inequality in the case of symmetric random walks. We extend this result in the centered case in the next proposition and then show in an example how this can be used to prove that some random walks have vanishing speed.

PROPOSITION 3.11. *Assume that* (C1) *holds. If the entropy vanishes, then $L = 0$.*

PROOF. It is straightforward once we recall the Carne–Varopoulos bound from Theorem (2.8): for some constant $C$, we have
$$\mu^t(x) \le C e^{-d^2(id,x)/(Ct)}.$$

For any $\alpha > 0$, we then have

$$0 = h = \lim -\frac{1}{t}\mathbb{E}[\log \mu^t(X_t)] = \lim -\frac{1}{t}\sum_{x \in V} \log \mu^t(x)\mu^t(x)$$

$$\geq \lim -\frac{1}{t} \sum_{x; d(id,x) \geq \alpha t} \log \mu^t(x)\mu^t(x)$$

$$\geq \lim \frac{1}{t} \sum_{x; d(id,x) \geq \alpha t} \left(-\log C + \frac{d^2(id,x)}{Ct}\right)\mu^t(x)$$

$$\geq \frac{\alpha^2}{C} \lim \sum_{x; d(id,x) \geq \alpha t} \mu^t(x) = \frac{\alpha^2}{C}\mathbb{P}[X_t \geq \alpha t].$$

Therefore $\mathbb{P}[X_t \geq \alpha t] \to 0$ for any $\alpha > 0$ and thus $L = 0$. $\square$

EXAMPLE 3.12. We discuss the application of the last proposition in the case of $\mathbb{Z} \wr \mathbb{Z}$ using the same notation as in the proof of Proposition (3.3)(c).

Let $G$ be a finite generating sequence in $\mathbb{Z} \wr \mathbb{Z}$ satisfying condition (C1) or equivalently condition (C2). Since $G$ generates $\mathbb{Z} \wr \mathbb{Z}$, its image by the homomorphism $\varepsilon$ generates $\mathbb{Z}$. Therefore the random walk $\varepsilon(X_t)$ is recurrent. It then follows from the fact that $\mathbb{Z} \wr \mathbb{Z}$ is a semidirect product of a recurrent group and an Abelian group that the Poisson boundary is trivial (see [8], Theorem 3.1), and therefore $h = 0$ and therefore, applying our proposition, $L = 0$.

It should be noted that if we drop the assumption that $G$ generates, the situation becomes quite different. Choose, for instance, $G = \{g_1 = (+2, T_1\sigma_1), g_2 = (-2, \sigma_{-1})\}$. Then $G$ satisfies condition (C2) since $\varepsilon(g_1) + \varepsilon(g_2) = +2 - 2 = 0$ and $\sum_x H(g_1)(x) + H(g_2)(x) = \sum_x \sigma_1(x) + \sigma_{-1}(x) = 0$. Clearly, $G$ does not satisfy condition (C1). As a matter of fact, there is no way to write $id$ as a nonempty product of $g_1$ and $g_2$. Besides $L \neq 0$. Indeed, each multiplication by $g_1$ adds a "1" at an odd location in $\mathbb{Z}$ and each multiplication by $g_2$ adds a "$-1$" at an even location in $\mathbb{Z}$. Thus $\sum_{x \in \mathbb{Z}} H(X_t)(2x) = -\#\{s \leq t : X_{s-1}^{-1} X_s = g_1\}$ and similarly $\sum_{x \in \mathbb{Z}} H(X_t)(2x+1) = \#\{s \leq t : X_{s-1}^{-1} X_s = g_2\}$. So $\sum_{x \in \mathbb{Z}} H(X_t)(2x+1) - H(X_t)(2x) = t$ and $L > 0$.

**Acknowledgments.** Many thanks to Christophe Pittet and Laurent Saloff-Coste for their help, suggestions and interest in this work. I also thank the referees for their comments.

## REFERENCES

[1] ALEXOPOULOS, G. K. (2002). Random walks on discrete groups of polynomial volume growth. *Ann. Probab.* **30** 723–801. MR1905856

[2] BALDI, P., LOHOUÉ, N. and PEYRIÈRE, J. (1977). Sur la classification des groupes récurrents. *C. R. Acad. Sci. Paris Ser. A–B* **285** 1103–1104. MR0518008

[3] CARNE, T. K. (1985). A transmutation formula for Markov chains. *Bull. Sci. Math.* **109** 399–405. MR0837740

[4] CHEN, M. F. (1991). Comparison theorems for Green functions of Markov chains. *Chinese Ann. Math. Ser. A* **12** 237–242. MR1130250

[5] DERRIENNIC, Y. (1980). Quelques applications du théorème ergodique sous-additif. *Astérisque* **74** 183–201. MR0588163

[6] HEBISCH, W. and SALOFF-COSTE, L. (1993). Gaussian estimates for Markov chains and random walks on groups. *Ann. Probab.* **21** 673–709. MR1217561

[7] JIANG, D. Q., QIAN, M. and QIAN, M. P. (2004). *Mathematical Theory of Nonequilibrium Steady States. On the Frontier of Probability and Dynamical Systems. Lecture Notes in Math.* **1833**. Springer, Berlin. MR2034774

[8] KAIMANOVICH, V. A. (1991). Poisson boundaries of random walks on discrete solvable groups. In *Probability Measures on Groups X* (H. Heyer, ed.) 205–238. Plenum, New York. MR1178986

[9] KAIMANOVICH, V. A. (2001). Poisson boundary of discrete groups. Available at name.math.univ-rennes1.fr/vadim.kaimanovich/.

[10] KALPAZIDOU, S. L. (1995). *Cycle Representations of Markov Processes*. Springer, New York. MR1336140

[11] VAROPOULOS, N. TH. (1985). Long range estimates for Markov chains. *Bull. Sci. Math.* **109** 225–252. MR0822826

[12] VERSHIK, A. M. (2000). Dynamic theory of growth in groups: Entropy, boundaries, examples. *Russ. Math. Surveys* **55** 667–733. MR1786730

[13] WOESS, W. (2000). *Random Walks on Infinite Graphs and Groups.* Cambridge Univ. Press. MR1743100

CMI
39 RUE JOLIOT-CURIE
13013 MARSEILLE
FRANCE
E-MAIL: pierre.mathieu@cmi.univ-mrs.fr